\documentstyle[11pt]{article}
\newtheorem{theorem}{Theorem}[section]
\newtheorem{lemma}[theorem]{Lemma}

\newtheorem{defi}[theorem]{Definition}
\def\f{\noindent}
\def\h{\hskip4pt  $\Box$ \vspace{4mm}}
\def\p{\f {\bf Proof}\hskip10pt}

\begin{document}
\title{The $k$th Upper Bases of Primitive Non-powerful Signed
Digraphs\thanks{Research
supported by NNSF of China (No. 10571163) and NSF of Shanxi (No. 20041010).}\\
\author{Yanling Shao${}^a$, \ Jian Shen${}^b$\thanks{Corresponding author.
E-mail addresses: ylshao@nuc.edu.cn (Y. Shao), js48@txstate.edu (J. Shen),
ybgao@nuc.edu.cn (Y. Gao)}, \ Yubin Gao${}^a$ \\
\small ${}^a$Department of Mathematics, North University of China,
Taiyuan, Shanxi 030051, P.R. China\\
\small ${}^b$Department of Mathematics, Texas State University,
San Marcos, TX 78666, USA  }}
\date{}
\maketitle

\begin{abstract}
In this paper, we study the $k$th upper bases
of primitive non-powerful signed digraphs.
A bound on the $k$th upper bases
of all primitive non-powerful signed digraphs is obtained,
and the equality case of the bound is characterized.
We also show that there exists ``gap" in the $k$th upper base set of
primitive non-powerful signed digraphs.

\bigskip
\f {\it AMS classification:}  05C20, 05C50, 15A48

\f {\bf Keywords:} Base, Sign pattern, Digraph, Primitive digraph
\end{abstract}

\section{Introduction}
\hskip\parindent
A {\em sign pattern} ({\em matrix}) $A$ is a
matrix whose entries are from the set $\{1, -1, 0\}$.
In the computations of the powers of a square sign pattern,
an ambiguous sign may arise
when a positive sign is added to a negative sign.
So a new symbol $\#$ was introduced in
\cite{LAA299} to denote such an ambiguous sign.
$\Gamma=\{0,1,-1,\#\}$ is called the {\em generalized sign set}
and addition and multiplication involving the symbol $\#$ are defined as follows:
$$(-1)+1=1+(-1)=\#; \ \ a+\#=\#+a=\# \ (\hbox{ for all } a \in \Gamma);$$
$$0\cdot \#=\#\cdot 0=0; \ \ b\cdot \#=\#\cdot b=\# \ (\hbox{ for all } b \in \Gamma \backslash \{0\}).$$
Following \cite{LAA299}, we call a matrix with entries from the set $\Gamma$
a {\em generalized sign pattern (matrix)}.

For a sign pattern $A$, we use $|A|$ to denote the $(0,1)$-matrix obtained from
$A$ by replacing each nonzero entry by $1$. Clearly
$|A|$ completely determines the zero pattern of $A$.
A square nonnegative matrix $A$ is {\em primitive}
if some power $A^k>0$ (that is, $A^k$ is entrywise
positive). The least such $k$ is called the {\em primitive exponent} of $A$,
denoted by ${\rm exp}(A)$. A square sign pattern
$A$ is called primitive if $|A|$ is primitive,
and in this case we define
${\rm exp}(A)={\rm exp}(|A|)$.

A {\em signed digraph} $S$ is a digraph where each arc is
assigned a sign $1$ or $-1$.
A {\em walk} $W$ in a digraph is a sequence of arcs, $e_1,e_2,\dots , e_k$, such that the terminal
vertex of $e_i$ is the same as the initial vertex of $e_{i+1}$ for $i=1,\dots ,k-1$.
The number $k$ of arcs is
called the {\em length} of the walk $W$, denoted by $l(W)$.
The {\em sign} of the walk $W$ is defined to be $\prod_{i=1}^k{\rm sgn}(e_i)$,
denoted by ${\rm sgn}W$.
Two walks $W_1$ and $W_2$ in a signed digraph
are called {\em a pair of $SSSD$ walks},
if they have the same initial vertex, same terminal vertex, same
length, but different signs.

Let $A$ be a sign pattern of order $n$ and $A, A^2, A^3, \dots$
be the sequence of powers of $A$. Since there are only $4^{n^2}$ different
generalized sign patterns of order $n$, there exist repetitions in
the above sequence. Let $l$ be the least positive integer
such that $A^l=A^{l+p}$ for some $p$.
Then $l$ is called the {\em generalized base} (or simply base) of $A$,
denoted by $l(A)$.
A square sign pattern $A$ is called {\em powerful} if each power
of $A$ contains no $\#$ entry.


Let $A=(a_{ij})$ be a sign pattern of order $n$.
The associated digraph $D(A)$ of $A$
is defined to be the digraph with vertex set $V=\{1,2,\dots ,n\}$ and arc
set $E=\{(i,j)|a_{ij} \not =0\}$. The {\em associated signed digraph}
$S(A)$ of $A$ is obtained from $D(A)$ by
assigning the sign of each $a_{ij}$ to the arc $(i,j)$ in $D(A)$.
Conversely, if
$S$ is a signed digraph, the sign pattern
$A$ of $S$ can be defined similarly.
We say that $S$ is primitive if $A$ is primitive, and
$S$ is powerful if $A$ is powerful. By \cite{DM}, $S$ is powerful if and only if
$S$ contains no pairs of $SSSD$ walks.
Also we define $l(S)=l(A)$.

Let $S$ be a primitive non-powerful signed digraph
of order $n$ and $k$ be an integer
with $1 \leq k \leq n$. For a subset $X \subseteq V(S)$ with $|X|=k$,
we define the {\em base of $X$}, denoted by $l_S(X)$,
to be the least integer $p$ such
that for each vertex $v \in V(S)$ there exists a pair of $SSSD$ walks of
length $p$ from some vertex in $X$ to $v$.
The {\em $k$th upper base} of $S$ is defined to be
$L(S,k)=\max\{l_S(X)|X \subseteq V(S) \hbox{ and }|X|=k \}$.
Clearly $L(S,1)=l(S)$ and $L(S,1) \geq L(S,2) \geq
\dots \geq L(S,n)$.

In this paper, we study the $k$th upper bases
of primitive non-powerful signed digraphs.
A bound on the $k$th upper bases
of all primitive non-powerful signed digraphs is obtained,
and the equality case of the bound is characterized.
We also show that there exists ``gap" in the $k$th upper base set of
primitive non-powerful signed digraphs.

\section{Preliminaries }
\begin{lemma}(\cite{DM}) \ \label{lemma21}
Let $S$ be a primitive signed digraph. Then $S$ is non-powerful if and only if $S$
contains a pair of cycles $C_1$ and $C_2$ (say, with lengths $p_1$ and $p_2$,
respectively) satisfying one of the following two conditions:

{\em (1)} $p_1$ is odd, $p_2$ is even, and ${\rm sgn} C_2=-1$;

{\em (2)} Both $p_1$ and $p_2$ are odd and ${\rm sgn} C_1=-{\rm sgn} C_2$.
\end{lemma}

For convenience, we call a pair of cycles $C_1$ and $C_2$
satisfying $(1)$ or $(2)$ in Lemma \ref {lemma21}
a {\em distinguished cycle pair}. Suppose
$C_1$ and $C_2$ form a distinguished cycle pair of lengths $p_1$ and $p_2$,
respectively. Then
the closed walks $W_1=\frac{{\rm lcm}(p_1,p_2)}{p_1}C_1$
(walk around $C_1$ $\frac{{\rm lcm}(p_1,p_2)}{p_1}$ times)
and $W_2=\frac{{\rm lcm}(p_1,p_2)}{p_2}C_2$ have the same length
${\rm lcm}(p_1,p_2)$ but with different signs since
$({\rm sgn} C_1)^{{\rm lcm}(p_1,p_2)/p_1}
=-({\rm sgn} C_2)^{{\rm lcm}(p_1,p_2)/p_2}.$
Similarly the closed walks $p_2C_1$
and $p_1C_2$ have the same length
$p_1p_2$ but with different signs.

Let $a_1,\dots, a_k$ be positive integers with ${\rm gcd}(a_1,\dots, a_k)=1$.
Define the Frobenius set as follows:
$$S(a_1,\dots, a_k)=\{r_1a_1+\dots+r_ka_k \ | \ r_1,\dots, r_k \hbox { \ are nonnegative integers}\}.$$
It is well-known that $S(a_1,\dots, a_k)$
contains all sufficiently
large positive integers. The Frobenius number
$\phi(a_1,\dots, a_k)$ is the least
integer $\phi$ such that $m\in S(a_1,\dots, a_k)$ for all integers $m \geq \phi$.
In particular, if $k=2$, it is known that
$\phi(a_1, a_2)=(a_1-1)(a_2-1)$.

\begin{defi}(\cite{JG14}) \ \label{defi22}
{\em Let $D$ be a primitive digraph of order $n$ and $k$ be an integer
with $1 \leq k \leq n$. For a subset $X \subseteq V(D)$ with $|X|=k$, we define
the} exponent of $X$, {\em denoted by ${\rm exp}_D(X)$,
 to be the smallest integer $p$ such that for each
vertex $v$ in $V(D)$ there exists a walk of length $p$ from some
vertex in $X$ to $v$.}
\end{defi}

\begin{defi}(\cite{JG14}) \ \label{defi23}
{\em Let $D$ be a primitive digraph of order $n$ and $1 \leq k \leq n$.
$$F(D,k)=\max\{{\rm exp}_D(X)|X \subseteq V(D) \hbox{ and } |X| = k\}$$
is called the} $k$th upper multiexponent {\em of $D$.}
\end{defi}

Thus for any $X \subseteq V(D)$ with $|X|=k$ and any $v$ in $V(D)$,
there exists a walk of length $l$ from some
vertex in $X$ to $v$ for all $l\geq F(D,k)$.


\setlength{\unitlength}{0.9mm}
\begin{picture}(180,49)
\put(45,40){\makebox(0,0){$\bullet$}}
\put(25,40){\makebox(0,0){$\bullet$}}
\put(25,10){\makebox(0,0){$\bullet$}}
\put(45,10){\makebox(0,0){$\bullet$}}
\put(55,25){\makebox(0,0){$\bullet$}}
\put(15,25){\makebox(0,0){$\bullet$}}
\put(25,40){\line(1,0){20}}
\put(45,40){\vector(-1,0){16}}
\put(45,10){\line(2,3){10}}
\put(45,10){\vector(2,3){8}}
\put(25,10){\line(-2,3){10}}
\put(15,25){\vector(2,-3){8}}
\put(45,40){\line(2,-3){10}}
\put(55,25){\vector(-2,3){8}}
\put(25,40){\line(-2,-3){10}}
\put(25,40){\vector(-2,-3){8}}
\put(55,25){\line(-2,1){30}}
\put(55,25){\vector(-2,1){24}}
\put(35,10){\makebox(0,0){$\cdots \cdots$}}
\put(48,42){\makebox(0,0){\small $n$}}
\put(22,42){\makebox(0,0){\small $1$}}
\put(12,25){\makebox(0,0){\small $2$}}
\put(63,25){\makebox(0,0){\small $n-1$}}
\put(53,10){\makebox(0,0){\small $n-2$}}
\put(21,10){\makebox(0,0){\small $3$}}
\put(38,1){\makebox(0,0){Fig. 2.1  \ \ The digraph $D_1$ }}
\put(125,40){\makebox(0,0){$\bullet$}}
\put(105,40){\makebox(0,0){$\bullet$}}
\put(105,10){\makebox(0,0){$\bullet$}}
\put(125,10){\makebox(0,0){$\bullet$}}
\put(135,25){\makebox(0,0){$\bullet$}}
\put(95,25){\makebox(0,0){$\bullet$}}
\put(105,40){\line(1,0){20}}
\put(125,40){\vector(-1,0){16}}
\put(125,10){\line(2,3){10}}
\put(125,10){\vector(2,3){8}}
\put(105,10){\line(-2,3){10}}
\put(95,25){\vector(2,-3){8}}
\put(125,40){\line(2,-3){10}}
\put(135,25){\vector(-2,3){8}}
\put(105,40){\line(-2,-3){10}}
\put(105,40){\vector(-2,-3){8}}
\put(125,40){\line(-2,-1){30}}
\put(125,40){\vector(-2,-1){24}}
\put(135,25){\line(-2,1){30}}
\put(135,25){\vector(-2,1){24}}
\put(115,10){\makebox(0,0){$\cdots \cdots$}}
\put(128,42){\makebox(0,0){\small $n$}}
\put(102,42){\makebox(0,0){\small $1$}}
\put(92,25){\makebox(0,0){\small $2$}}
\put(143,25){\makebox(0,0){\small $n-1$}}
\put(133,10){\makebox(0,0){\small $n-2$}}
\put(101,10){\makebox(0,0){\small $3$}}
\put(118,1){\makebox(0,0){Fig. 2.2  \ \ The digraph $D_2$ }}
\end{picture}

In the literature, many people defined $F(D,n)=1$
to emphasize that $F(D,n)$ must be a positive integer.
Here we define $F(D,n)=0$ since, for each vertex $u$,
there is a degenerated path of length $0$ from $u$ to itself.
The main purpose of this adjusted definition is to
include the trivial case $k=n$ into the following three lemmas
which were originally stated for $1 \leq k \leq n-1$.
Our results in Section $3$ are for all $k$, $1 \leq k \leq n$.

\begin{lemma}(\cite{JG18}) \ \label{lemma24}
Let $1 \leq k \leq n$ and $D$ be a primitive digraph of order $n$.
Then $F(D,k) \leq (n-k)(n-1)+1$.
\end{lemma}

\begin{lemma}(\cite{JG18}) \ \label{lemma26}
Let $1 \leq k \leq n$ and $D$ be a primitive digraph of order $n$ and $s$ be the length of a
shortest cycle of $D$. Then
$F(D,k) \leq (n-k-1)s+n$.
\end{lemma}

\begin{lemma}(\cite{GC14}) \ \label{lemma25}
Let $1 \leq k \leq n$ and $D_2$ be the primitive digraph defined in Fig. 2.2.
Then $F(D_2,k)=(n-k)(n-1)$.
\end{lemma}

\section{Main Results }
\begin{lemma}\ \label{lemma33}
Let $1 \leq k \leq n$ and $D_1$ be the primitive digraph defined in Fig. 2.1.
Let $X$ be a subset of $V(D_1)$ with $|X|=k$.
Then, for any integer $l \geq (n-k)(n-1)$ and any vertex $i$ with $1 \leq i \leq n-1$,
there exists a walk of length $l$ from some $x_0 \in X$ to $i$ in $D_1$.
\end{lemma}

\p Denote $p=l-(n-k)(n-1)$.
Then there exists
a walk of length $p$ from some vertex $j$, $1 \leq j \leq n-1$, to the vertex $i$.
Let $A$ be the adjacency matrix of $D_1$. Then
there exists a loop at vertex $j$
in $D(A^{n-1})$, where $D(A^{n-1})$ is the associated digraph of $A^{n-1}$.
For any $X\subseteq V(D_1)$ with $|X|=k$,
there exists a walk of length
$n-k$ from some $x_0\in X$ to $j$ in $D(A^{n-1})$.
So there exists a walk of length $(n-k)(n-1)$ from $x_0$ to $j$
in $D_1$, and thus this walk can be extended to a walk of length
$l$ from $x_0$ to $i$ in $D_1$. \h

\begin{lemma}\ \label{lemma31}
Let $S$ be a primitive non-powerful signed digraph of order $n$ with $D$
as its underlying digraph.
Let $W_1$ and $W_2$ be a pair of $SSSD$ walks of length
$r$ from a vertex $u$ to a vertex $v$. Then, for $1 \leq k \leq n$,
$$L(S,k) \leq F(D,k)+d(D)+r,$$
where $d(D)$ is the diameter of $D$.
\end{lemma}

\p
Let $X \subseteq V(S)$ with $|X|=k$ and $y$ be any vertex of $S$.
Let $P$ be a shortest path in $S$ from $v$ to $y$.
Then $l(P) \leq d(D)$.
Denote $q=F(D, k)+d(D)-l(P)$. Then $q \geq F(D, k)$.
By Definitions \ref{defi22} and \ref{defi23},
there exists a walk $Q$ of length
$q$ from some $x_0 \in X$ to $u$.
Therefore $Q+W_1+P$ and $Q+W_2+P$ form a pair of $SSSD$ walks
of length $q+r+l(P)=F(D,k)+d(D)+r$
from $x_0$ to $y$.
So $l_S(X) \leq F(D,k)+d(D)+r$, from which Lemma \ref{lemma31}
follows.
\h

\begin{lemma}\ \label{lemma32}
Let $S$ be a primitive non-powerful signed digraph of order $n$ with $D$
as its underlying digraph. Let $V_1\subseteq V(S)$.
If for each vertex $u \in V_1$, there is
a pair of $SSSD$ walks of length $r$ from $u$ to $u$, then
for $1 \leq k \leq n$,
$$L(S,k) \leq F(D,k)+r+n-|V_1|.$$
\end{lemma}

\p
By Lemma \ref{lemma31}, we may suppose $V_1\not = \phi$.
Let $X \subseteq V(S)$ with $|X|=k$ and $y$ be any vertex of $S$.
Then there is some vertex $v \in V_1$ such that
the shortest path $P$ from $v$ to $y$ has length at most
$n-|V_1|$. By Definitions \ref{defi22} and \ref{defi23},
there
exists a walk $Q$ of length $F(D,k)$ from some $x_0\in X$ to $v$.
Let $W_1$ and $W_2$
be a pair of $SSSD$ walks of length $r$ from $v$ to $v$.
Then $Q+W_1+P$ and $Q+W_2+P$ form a pair of $SSSD$ walks of
length $F(D, k)+r+n-|V_1|$ from $x_0$ to $y$. So $L(S,k) \leq F(D,k)+r+n-|V_1|$.
\h

\begin{theorem} \ \label{th33}
Let $S_1$ be a primitive non-powerful signed digraph of order $n$ with $D_1$ (Fig. 2.1)
as its underlying digraph. Then for $1 \leq k \leq n$,
$$L(S_1,k)=(2n-k)(n-1)+1.$$
\end{theorem}

\p
First we show that there is a pair of $SSSD$ walks of length $(n-1)^2+1$
from vertex $n-1$ to $1$. Let $Q_1$ and $Q_2$ be the paths of lengths $1$ and $2$
from the vertex $n-1$ to vertex $1$, respectively.
Let $C_n$ and $C_{n-1}$ be the cycles of lengths
$n$ and $n-1$, respectively. Take $W_1=Q_1+(n-1)C_{n-1}$, and $W_2=Q_2+(n-2)C_n$.
Let $P_1$ be the unique path from vertex $1$ to $n-1$. Then
$W_1+P_1=nC_{n-1}$, and $W_2+P_1=(n-1)C_n$.
Because $S_1$ is primitive non-powerful and $S_1$ has only two cycles
of lengths $n$ and $n-1$ respectively, by Lemma \ref{lemma21},
${\rm sgn}(nC_{n-1})=-{\rm sgn}((n-1)C_n)$. Thus
${\rm sgn}(W_1)=-{\rm sgn}(W_2)$ and $W_1$ and $W_2$ form a
pair of $SSSD$ walks of length $(n-1)^2+1$.

Let $X \subseteq V(S_1)$ with $|X|=k$ and $y$ be any vertex of $S_1$.
Let $P_2$ be the shortest path from $1$ to $y$.
Then $l(P_2) \leq n-1$.
Denote $q=(n-k)(n-1)+n-1-l(P_2)$. Then $q \geq (n-k)(n-1)$.
By Lemma \ref{lemma33},
there exists a walk $Q$ of length
$q$ from some $x_0 \in X$ to $n-1$.
Therefore $Q+W_1+P_2$ and $Q+W_2+P_2$ form a pair of $SSSD$ walks of length
$(2n-k)(n-1)+1$ from $x_0$ to $y$.
So $L(S_1,k) \leq (2n-k)(n-1)+1$.

On the other hand, let $X_0=\{n,n-1, \dots, n-k+1\}$ such that $|X_0|=k$.
We will show that there is no pair of $SSSD$ walks of length
$(2n-k)(n-1)$ from any vertex of $X_0$ to $n$.
Suppose $W_1$ and $W_2$ are two walks of length $(2n-k)(n-1)$ from the vertex
$n-m$ to $n$, where $0 \leq m \leq k-1$. Then each $W_i$ is a union of
the unique path $P$ from $n-m$ to $n$ of length $m$ and a number of
cycles of lengths $n-1$ and $n$. Thus
there exist nonnegative integers $a_i, b_i \ (i=1,2)$ with $a_2\geq a_1$ such that
$$(2n-k)(n-1)=a_in+b_i(n-1)+m, \ \ \ (i=1,2).$$
This implies that $\phi (n-1,n)-1=(n-1)(n-2)-1
=(a_2+m-n)n+(b_2+k-m-1)(n-1)$.
Also we have $(a_2-a_1)n=(b_1-b_2)(n-1)$ and thus $(n-1)|(a_2-a_1)$.
If $a_2 \not = a_1$, then
$a_2+m-n \geq a_1+n-1+m-n=a_1+m-1\geq 0$,
where the last inequality holds since $a_1 \geq 1$ (that is,
$W_1$ contains at least one cycle of length $n$) if $m=0$.
Recall that $m\leq k-1$. Then $b_2+k-m-1 \geq 0$,
which contradicts that $\phi (n-1,n)-1$
can not be expressed as a non-negative integral linear
combination of $n-1,n$. So $a_1=a_2$
and thus ${\rm sgn}(W_1)={\rm sgn}(W_2)$. This shows that there is
no pair of $SSSD$ walks of length $(2n-k)(n-1)$
from vertex $n-m$ to $n$ in $S_1$ for $m=0,1, \dots , k-1$.
So $L(S_1, k)\geq (2n-k)(n-1)+1$, from which Theorem \ref{th33} follows.
\h

\begin{theorem} \ \label{th34}
Let $S_2$ be a primitive non-powerful signed digraph of order $n$ with $D_2$ (Fig. 2.2)
as its underlying digraph.

{\em (1)} \ If the (only) two cycles of length $n-1$ of $S_2$
have different signs, then for $1 \leq k \leq n$,
$$L(S_2, k) \leq (n-k+1)(n-1)+2.$$

{\em (2)} \ If the (only) two cycles of length $n-1$ of $S_2$
have the same sign, then for $1 \leq k \leq n$,
$$L(S_2,k)=(2n-k)(n-1).$$
\end{theorem}

\p
(1) The (only) two cycles of length $n-1$ of $S_2$ have different signs.
Let $Q_1=(n-1,1)+(1,2)$ and $Q_2=(n-1,n)+(n,2)$ be the two paths of length
$2$ from the vertex $n-1$ to the vertex $2$.
Then ${\rm sgn}Q_1=-{\rm sgn}Q_2$. By Lemmas \ref{lemma31} and \ref{lemma25},
$L(S_2,k) \leq F(D_2,k)+d(D_2)+r\leq (n-k)(n-1)+(n-1)+2
=(n-k+1)(n-1)+2$.

(2) The two cycles of length $n-1$ in $S_2$ have the same sign.
Then each cycle of length $n-1$ and the unique cycle of length $n$
form a distinguished cycle pair. By Lemma \ref{lemma21},
there is a pair of $SSSD$ walks
of length $n(n-1)$ from each $i$ to $i$ $(1 \leq i \leq n)$.
By Lemma \ref{lemma32},
$L(S_2, k) \leq F(D_2,k)+n(n-1) \leq (n-k)(n-1)+n(n-1)=(2n-k)(n-1)$.

On the other hand, let $X_0=\{1,n,n-1, \dots ,n-k+2\}$ such that $|X_0|=k$.
(If $k=1$, then $X_0=\{1\}$.) 
We will show that there is no pair of $SSSD$ walks of length $(2n-k)(n-1)-1$
from any vertex of $X_0$ to $n$.
Suppose $W_1$ and $W_2$ are two walks of length $(2n-k)(n-1)-1$ from vertex
$n-m$ to $n$, where either $m=n-1$ or $0\leq m \leq k-2$.
Then each $W_i$ is a union of the unique path $P$ of length $m$ from $n-m$ to $n$ and
a number of cycles of lengths $n-1$ and $n$. Thus there exist
nonnegative integers $a_i, b_i  \ (i=1,2)$ with $b_1\geq b_2$ such that
$$(2n-k)(n-1)-1=a_in+b_i(n-1)+m, \  (i=1,2).$$
This implies that

\hspace{15mm} $\phi (n-1,n)-1=(n-1)(n-2)-1$

\hspace{2.58cm} $$= \left\{\begin{array}{lcl}
a_1n+(b_1+k-n-1)(n-1) &&\hbox{if } m=n-1;\\
(a_1+m)n+(b_1+k-m-2-n)(n-1) &&\hbox{if } 0 \leq m \leq k-2.
\end{array}\right.$$
Also we have $(a_2-a_1)n=(b_1-b_2)(n-1)$ and thus $n|(b_1-b_2)$.
If $b_1\not =b_2$, then $b_1 \geq n$ and thus
$b_1+k-n-1\geq 0$ and
$b_1+k-m-2-n \geq 0$ when $0 \leq m \leq k-2$.
This contradicts that $\phi (n-1,n)-1$
can not be expressed as a non-negative integral
linear combination of $n-1,n$. So $b_1=b_2$
and thus
${\rm sgn}(W_1)={\rm sgn}(W_2)$.
This shows that there is
no pair of $SSSD$ walks of length $(2n-k)(n-1)-1$
from any vertex of $X_0$ to $n$.
So $L(S_2, k)\geq (2n-k)(n-1)$, from which Theorem \ref{th34} follows.\h

\begin{theorem} \ \label{th35}
Let $S$ be a primitive non-powerful signed digraph of order $n \geq 6$ with $D$
as its underlying digraph, where $D$ is not isomorphic to $D_1$ or $D_2$ (Fig. 2.1,2.2).
Then for $1 \leq k \leq n$,
$$L(S, k) \leq (2n-k)(n-1)-3.$$
\end{theorem}

\p
Let $C_s$ be a shortest cycle of length $s$ in $S$.
Then $s \leq n-2$
since $D$ is not isomorphic to $D_1$ or $D_2$. By Lemma \ref {lemma26},
$F(D,k) \leq (n-k-1)s+n \leq (n-k-1)(n-2)+n$.
Because $S$ is primitive non-powerful, by Lemma \ref{lemma21},
there is a distinguished cycle pair $C_1$ and $C_2$
(of lengths, say, $p_1$ and $p_2$ respectively) in $S$. Then
$p_1C_2$ and $p_2C_1$ have different signs.
We may assume $p_1\geq p_2$.

Case 1. \ $\{p_1, p_2\}\not =\{n-1, n\}$.
Let $V_1=V(C_1)\cap V(C_2)$. If $V_1=\phi$,
then there is some
$i \in V(C_1)$ and $j \in V(C_2)$ such that
$d(i,j)\leq n-p_1-p_2+1$. Then there
is a pair of $SSSD$ walks of length $d(i,j)+{\rm lcm}(p_1,p_2)$
from $i$ to $j$. By Lemma \ref{lemma31},
$L(S,k) \leq F(D, k)+d(D)+d(i,j)+{\rm lcm}(p_1,p_2)
\leq F(D, k)+2n-p_1-p_2+{\rm lcm}(p_1,p_2)$.
On the other hand, if $V_1\not =\phi$,
then $|V_1|\geq |V(C_1)|+|V(C_2)|-n=p_1+p_2-n$.
Since for each $u\in V_1$, there is a pair of $SSSD$ walks of
length ${\rm lcm}(p_1,p_2)$ from $u$ to $u$, by Lemma \ref{lemma32},
$L(S,k) \leq F(D, k)+{\rm lcm}(p_1,p_2)+n-|V_1|
\leq F(D, k)+2n-p_1-p_2+{\rm lcm}(p_1,p_2)$.
Thus for all possible $V_1$, we always have

\hspace{0mm} $L(S,k) \leq F(D, k)+2n-p_1-p_2+{\rm lcm}(p_1,p_2)$

$$\leq \left\{\begin{array}{lcl}
(n-k-1)(n-2)+n+2n+(n-1)(n-3)-1 &&\hbox{if } k\leq n-1\\
0+2n+(n-1)(n-3)-1 &&\hbox{if } k=n
\end{array}\right.$$

\hspace{1.3cm} $ \leq (2n-k)(n-1)-3.$

Case 2. \ $p_1=n$ and $p_2=n-1$.

Subcase 2.1. \ $s$ is odd. Then either
$C_1$ and $C_s$ (when $n$ is even) or $C_2$ and
$C_s$ (when $n$ is odd) also  form a distinguished cycle pair. By
Case 1, $L(S, k) \leq (2n-k)(n-1)-3$.

Subcase 2.2. \ $s$ is even and $\hbox{sgn}C_s=-1$. Then either $C_1$
and $C_s$ (when $n$ is odd) or $C_2$ and $C_s$ (when $n$ is even) also
form a distinguished cycle pair.
By Case 1, $L(S, k) \leq (2n-k)(n-1)-3$.

Subcase 2.3. \ $s$ is even and $\hbox{sgn}C_s=1$.
Then $|V(C_1)\cap V(C_2)\cap V(C_s)|=|V(C_2)\cap V(C_s)|
\geq p_2+s-n=s-1\geq 1$.
Let $j\in V(C_1)\cap V(C_2)\cap V(C_s)$. Then
$(n-s-1)C_1+C_s$ and $(n-s)C_2$
form a pair of $SSSD$
walks of length $(n-s)(n-1)$ from $j$ to $j$.
By Lemma \ref{lemma31},

\hspace{6mm} $L(S,k) \leq F(D, k)+d(D)+(n-s)(n-1)$

$$\leq \left\{\begin{array}{lcl}
(n-k-1)s+n+n-1+(n-s)(n-1) &&\hbox{if } k\leq n-1\\
0+n-1+(n-s)(n-1) &&\hbox{if } k=n
\end{array}\right.$$

\hspace{1.9cm} $ \leq (2n-k)(n-1)-3.$\h

By Theorems \ref{th33}, \ref{th34} and \ref{th35}, we obtain
our main theorem.

\begin{theorem} \ \label{th36}
Let $1 \leq k \leq n$ and $S$ be a primitive non-powerful
signed digraph of order $n \geq 6$. Then,

{\em (1)} \ $L(S, k) \leq (2n-k)(n-1)+1$.

{\em (2)} \ $L(S, k)=(2n-k)(n-1)+1$ if and only if the underlying digraph
of $S$ is isomorphic to $D_1$.

{\em (3)} \ $L(S, k)=(2n-k)(n-1)$ if and only if the underlying digraph
of $S$ is isomorphic to $D_2$ and the two cycles of length $n-1$ have the same sign.

{\em (4)} \ For each integer $t$ with $(2n-k)(n-1)-3<t<(2n-k)(n-1)$,
there is no primitive non-powerful signed digraph $S$ of order $n$ with
$L(S,k)=t$.
\end{theorem}

\end {document}